\documentclass[12pt]{amsart}

\usepackage{amssymb}
\usepackage{amsmath}
\evensidemargin\oddsidemargin

\begin{document}

\setcounter{page}{1}

\newtheorem{REM}{Remark\!\!}
\newtheorem{REMS}{Remarks\!\!}
\newtheorem{LEMO}{Lemma 1\!\!}
\newtheorem{LEMU}{Lemma 2\!\!}
\newtheorem{THEO}{Theorem\!\!}
\newtheorem{DEFI}{Definition\!\!}
\renewcommand{\theTHEO}{}
\renewcommand{\theDEFI}{}
\renewcommand{\theLEMO}{}
\renewcommand{\theLEMU}{}
\renewcommand{\theREM}{}
\renewcommand{\theREMS}{}

\newcommand{\eqnsection}{
\renewcommand{\theequation}{\thesection.\arabic{equation}}
    \makeatletter
    \csname  @addtoreset\endcsname{equation}{section}
    \makeatother}
\eqnsection

\def\a{\alpha}
\def\b{\beta}
\def\ba{b_\a}
\def\CC{{\mathbb{C}}} 
\def\Eac{{\mathcal E}}
\def\EE{{\mathbb{E}}} 
\def\elaw{\stackrel{d}{=}}
\def\eps{\varepsilon}
\def\fa{f_\a}
\def\ga{g_\a}
\def\hS{{\hat S}}
\def\hT{{\hat T}}
\def\hX{{\hat X}}
\def\ii{{\rm i}}
\def\lbd{\lambda}
\def\lacc{\left\{}
\def\lcr{\left[}
\def\lpa{\left(}
\def\lva{\left|}
\def\NN{{\mathbb{N}}} 
\def\pb{{\mathbb{P}}}
\def\rl{{\mathbb{R}}}
\def\racc{\right\}}
\def\rcr{\right]}
\def\rpa{\right)}
\def\Un{{\bf 1}}
\def\Ya{Y_\a}
\def\Za{Z_\a}

\def\d{\, \mathrm{d}}
\def\qed{\hfill$\square$}
\def\elaw{\stackrel{d}{=}}

\newcommand{\fin}{\vspace{-0.4cm}
                  \begin{flushright}
                  \mbox{$\Box$}
                  \end{flushright}
                  \noindent}
                  
\title[Multiplicative unimodality for stable laws]
      {Multiplicative strong unimodality for positive stable laws}

\author[Thomas Simon]{Thomas Simon}

\address{Laboratoire Paul Painlev\'e, U. F. R. de Math\'ematiques, Universit\'e de Lille 1, F-59655 Villeneuve d'Ascq Cedex. {\em Email} : {\tt simon@math.univ-lille1.fr}}



\begin{abstract} It is known that real Non-Gaussian stable distributions
  are unimodal, not additive strongly unimodal, and multiplicative
  strongly unimodal in the symmetric case. By a theorem of
  Cuculescu-Theodorescu \cite{CT}, the only remaining relevant situation for the
  multiplicative strong unimodality of stable laws is the one-sided. In this paper, we show that positive $\a-$stable laws are multiplicative strongly unimodal iff $\a\le 1/2.$
\end{abstract}

\maketitle

\section{The MSU property and stable laws}
 
A real random variable $X$ is said to be unimodal (or quasi-concave) if there exists $a\in\rl$ such that the functions $\pb[X\le x]$ and $\pb[X>x]$ are convex respectively in $(-\infty,a)$ and $(a,+\infty).$ If $X$  is absolutely continuous, this means that its density increases on $(-\infty,a]$ and decreases on $[a,+\infty).$ The number $a$ is called a mode of $X$, and might not be unique. A well-known example due to Chung shows that unimodality is not stable under convolution, and for this reason the notion of strong unimodality had been introduced in \cite{I}: a real variable $X$ is said to be strongly unimodal if the independent sum $X + Y$ is unimodal for any unimodal variable $Y$ (in particular $X$ itself is unimodal, choosing $Y$ degenerated at zero). In \cite{I} Ibragimov also obtained the celebrated criterion that $X$ is strongly unimodal iff it is absolutely continuous with a log-concave density.

Proving unimodality or strong unimodality properties is simple for variables
with given densities, but the problem might turn out complicated when these
densities are not explicit. In this paper we will deal with real (strictly)
stable variables, where very few closed formulae (given e.g. in \cite{Z}
p. 66) are known for the densities. A classical theorem of Yamazato shows that
they are all unimodal with a unique mode - see Lemma 1 in \cite{IC} for the 
previously shown one-sided case and Theorem 53.1 in \cite{S} for the general
result, but except in the Gaussian situation it is easy to see that stable laws are not strongly unimodal, because their heavy tails prevent the densities from being everywhere log-concave - see Remark 52.8 in \cite{S}.

Having infinitely divisible distributions, stable variables appear naturally in additive identities, a framework where they hence may not preserve unimodality. Stable variables also occur in multiplicative factorizations as a by-product of the so-called $M$-scheme (or $M$-infinite divisibility),  a feature which has been studied extensively by Zolotarev - see Chapter 3 in \cite{Z} and also \cite{P} for the one-sided case. Another concrete example of a multiplicative identity involving a stable law is the following. Suppose that $X, Y$ are two positive variables whose Laplace transforms are connected through the identity $\EE [e^{-\lbd X}] \, = \, \EE [e^{-\lbd^\a Y}]$
for some $\a \in (0,1)$ (up to a reformulation one could also take $\a >1$). Then one has  
$$ X\; \elaw \; \Za\;\times\; Y^{1/\a}$$
where $\Za$ is an independent positive $\a-$stable variable. In such identities, one can be interested in understanding whether the factorization through $\Za$ modify or not some basic distributional properties. 

From the point of view of unimodality, it is therefore natural to ask whether
stable variables are multiplicative strongly unimodal, in other words, whether
their independent product with any unimodal variable remains unimodal, or
not. A quick answer can be given in the symmetric case  because the mode of a
symmetric stable variable is obviously zero: Khintchine's theorem entails then
that its product with {\em any} independent variable will be unimodal with
mode at zero - see Proposition 3.6 in \cite{CT} for details. However, in the
non-symmetric case the mode of a stable variable is never zero (this is
obvious in the drifted Cauchy case $\a = 1$ and we refer to \cite{Z} p. 140
for the case $\a\neq 1$), so that such a simple argument cannot be
applied. The following criterion established in \cite{CT} Theorem 3.7, is a multiplicative counterpart to Ibragimov's theorem:

\begin{THEO}[Cuculescu-Theodorescu] Let $X$ be a unimodal random variable such
  that 0 is not a mode of $X$. Then $X$ is multiplicative strongly unimodal if and only if it is one-sided and absolutely continuous, with a density $f^{}_{X}$ having the property that
\begin{equation}
\label{CCTD}
t\;\mapsto\; f^{}_{X}(e^t) \quad\mbox{is log-concave in $\rl$} 
\end{equation}
when $X$ is positive (resp. $t\,\mapsto\, f^{}_{X}(-e^t)$ is log-concave in $\rl$
when $X$ is negative). 
\end{THEO}
With a slight abuse of notation, in the following we will say that a positive
random variable is {\em MSU\,} if and only if (\ref{CCTD}) holds. Cases of
multiplicative strongly unimodal, positive  variables with mode at zero and
such that (\ref{CCTD}) does not hold, are hence excluded in this
definition. But these cases are particular, and by the above remark no more relevant to the content of this paper which deals with stable laws. Besides, a change of variable and Ibragimov's theorem entail the following useful characterization for positive variables:
\begin{equation}
\label{Ibr}
\mbox{$X$ is MSU}\;\;\Longleftrightarrow\;\; \mbox{$\log X$ is strongly unimodal.}
\end{equation}
In particular the MSU property is stable by inversion and also, from
Pr\'ekopa's theorem, by independent multiplication. Another important feature
coming from (\ref{Ibr}) is that the MSU property remains unchanged under
rescaling and power transformations, which also comes from the obvious
analytical fact that (\ref{CCTD}) holds iff $t \mapsto K_1 e^{a_1 t} f^{}_X (K_2 e^{a_2 t})$ is log-concave for some $a_1\in\rl$ and $a_2, K_1, K_2 \in\rl^*.$ Notice however that the MSU property is barely connected to the strong unimodality of $X$ itself (several examples of this difference are given in \cite{CT}).
 
In this paper we are interested in the MSU property for positive $\a-$stable laws. For every $\a \in\, ]0,1[,$ consider $\fa$ the positive $\a-$stable density and $\Za$ the corresponding random variable, normalized such that
\begin{equation}
\label{Laplace}
\int_0^\infty e^{-\lbd t} \fa(t) dt\; =\; \EE\lcr e^{-\lbd \Za}\rcr \; =\; e^{-\lbd^\a}, \qquad \lbd \ge 0.
\end{equation}
Before studying the MSU property for $\Za$, in view of (\ref{Ibr}) one must first ask if $\log \Za$ is simply unimodal. A positive answer for all $\a \in\, ]0,1[$ had been given by Kanter - see Theorem 4.1 in \cite{K2}, who also deduced from \cite{IC} the decomposition
$$\log \Za\; =\; \a^{-1}\log \ba(U)\; +\; (\a-1)/\a\log L$$
where $L$ is a standard exponential variable, $U$ an independent uniform
variable over $[0,\pi],$ and $b_\a(u) = (\sin(\a u)/\sin (u))^\a(\sin((1-\a)
u)/\sin(u))^{1-\a}.$ The random variable $\log L$ is easily seen to be
strongly unimodal, but $\log \ba(U)$ is not - at least for $\a = 1/2,$ see the
Remark before Theorem 4.1 in \cite{K}. This leaves the question of the MSU
property for positive stable laws unanswered, and our result aims at filling this gap:

\begin{THEO} The variable $\Za$ is MSU if and only if $\a \le 1/2.$
\end{THEO}

To conclude this introduction, let us give two further reformulations of (\ref{CCTD}) in the positive stable case. The first one lies at the core of our proof, whereas the second one is probably nothing but a mere curiosity. Since $\fa$ is smooth, differentiating twice the logarithm entails that (\ref{CCTD}) is equivalent to the inequality
\begin{equation}
\label{LCE}
(x^2 f_\a''(x) + xf_\a'(x))f_\a(x)\; \le \; x^2(f_\a'(x))^2, \quad x\ge 0.
\end{equation}
If $m_\a$ stands now for the mode of $\Za$ and if $x_\a = \inf \{ x > m_\a, \;
\fa''(x) = 0\},$ then we know from (53.13) in \cite{S} - an identity which is
proved there for a certain class of positive self-decomposable distributions, but it
can also be obtained for positive stable laws after taking the weak limit -
and the discussion thereafter that (\ref{LCE}) is true for any $x\in [0,
x_\a].$ Hence, the MSU property amounts to the fact that it remains true for
all $x > x_\a$. From (53.13) in \cite{S} we also know that (\ref{LCE}) is
equivalent to the positivity everywhere of the function
$$x\;\mapsto\; \int_0^x (\fa'(x-y)\fa(x) - \fa(x-y)\fa'(x))y^{-\a}dy,$$
but this criterion is not very tractable because of the long memory
involved in the integral.
Thanks to the Humbert-Pollard representation for $\fa$ which will be recalled
at the beginning of Section 3, we finally mention that (\ref{LCE}) is equivalent to
$$\lpa\sum_{n\ge 1} (1+\a n)^2 \frac{(-1)^n  x^{-(1 +\a n)}}{\Gamma (-n\a)n!}   \rpa\lpa\sum_{n\ge 1} \frac{(-1)^n x^{-(1 +\a n)}}{\Gamma (-n\a)n!} \rpa\; \le\; \lpa\sum_{n\ge 1}  (1+\a n) \frac{(-1)^n x^{-(1 +\a n)}}{\Gamma (-n\a)n!} \rpa^2,$$
a strange inequality which would plainly hold in the opposite direction if the terms of the series had constant signs.

\section{Some particular cases} 

In this section we depict some situations where the MSU property can be proved
or disproved directly, thanks to more or less explicit representations for the
density $\fa$ or the variable $\Za.$ First of all, in the case $\a = 1/2$ it
is readily seen from the known formula
$$f_{1/2}(x)\; =\; \frac{1}{2\sqrt{\pi x^3}}e^{-1/4x}$$
that $Z_{1/2}$ is MSU. When $\alpha = 1/3$, Formula (2.8.31) in \cite{Z} yields 
$$f_{1/3}(x)\; =\; \frac{1}{3\pi x^{3/2}} K_{1/3}\lpa \frac{2}{3\sqrt{3}} x^{-1/2}\rpa$$
where $K_{1/3}$ is the Macdonald function of order 1/3, so that (\ref{CCTD})
amounts to show that $t\mapsto K_{1/3} (e^t)$ is log-concave. From (\ref{LCE})
and since $K_{1/3}$ is a solution to the modified Bessel equation
$$x^2 K_{1/3}'' \;+\; x K_{1/3}'\; =\; (x^2 + 1/9)K_{1/3}$$
over $\rl$, this is equivalent to $(x^2 + 1/9)K_{1/3}^2(x)\,\le \, x^2(
K_{1/3}'(x))^2$ for every $x\ge 0.$ Because $K_{1/3}'(x) < 0$ for every $x \ge 0,$ the latter is now an immediate consequence of a Turan-type inequality for modified Bessel functions recently established in \cite{Ba} - see (2.2)
therein, so that $Z_{1/3}$ is MSU.

With a probabilistic and more concise argument, one can also show that Property (\ref{CCTD}) holds for every 
$\alpha = 1/p, \, p \ge 2.$ A classical representation originally due to
E.~J.~Williams - see Section 2 in \cite{W} - shows indeed that after rescaling
$Z_{1/p}^{-1}$ is an independent product of $p$ Gamma variables:
$$Z_{1/p}^{-1}\; \elaw\; p^p\;Y_{1/p}\;\times\; Y_{2/p}\;\times\;\cdots\;\times\; Y_{(p-1)/p}$$
where each $Y_{k/p}$ is MSU since its density is $\Gamma(k/p)^{-1} y^{k/p -1}
e^{-y}\Un_{\{y > 0\}}.$ Hence, since the MSU property is stable by inversion and
independent multiplication, it follows that $Z_{1/p}$ is MSU for every $p\ge 2.$ 

As a matter of fact, a much stronger property is known to hold when $\alpha = 1/p$ for some $p \ge 2.$ In theses cases, it had namely been noticed in \cite{K} thanks to the classical computation of the fractional moments
\begin{equation}
\label{MellinS}
\EE[\Za^s]\; =\; \frac{\Gamma(1 - s/\a)}{\Gamma(1-s)}
\end{equation}
for every $s < \a,$ and the duplication formula for the Gamma function, that the kernel $(x,y)\mapsto f_\a(e^{x-y})$ is totally positive - see pp. 121-122 and 390 in \cite{K} as well as pp. 11-12 therein for the definition of total positivity. In particular, it is totally positive of order 2 (${\rm TP}_2$) which means precisely that  $x\mapsto f_\a(e^x)$ is a log-concave function - see e.g. 
Theorem 4.1.9 in \cite{K} Chap. 4. However, when $\a$ is not the reciprocal of an integer, the function 
$$s\; \mapsto\;\frac{\Gamma(1 -s)}{\Gamma(1-s/\a)}$$ 
is not an entire function of the type $\Eac_2^*$ - see e.g. \cite{K} p. 336 for a definition - and by Theorem 7.3.2 in \cite{K}, this
entails that $(x,y) \mapsto f_\a(e^{x-y})$ is no more a totally positive 
kernel. Karlin raised then the question whether or not it should be 
totally positive of some finite order - see \cite{K} p. 390, a problem which
seems as yet unadressed.

Let us finally show that $Z_{2/3}$ is not MSU, in other words that $(x,y) \mapsto
f_{2/3}(e^{x-y})$ is not a ${\rm TP}_2$ kernel. Formula (2.8.33) in \cite{Z}
(with a slight normalizing correction therein) yields first the expression
$$f_{2/3}(x)\; =\; \sqrt{\frac{3}{\pi x}}\, e^{-2/27x^2} W_{1/2, 1/6} \lpa
\frac{4}{27} x^{-2}\rpa$$
where $W_{1/2, 1/6}$ is a Whittaker function. Hence, from formul\ae \,(6.9.2) and 
(6.5.2) in \cite{E} we see that (\ref{CCTD}) is equivalent to the log-concavity of $t\mapsto g(e^t)$ with $g(x) = e^{-x} U_4(x)$ and the notation $U_\lbd(x) =
\Psi(1/6, \lbd/3, x)$ for all $\lbd >1/2,$ where
$$\Psi(a, c, x)\; =\; \frac{1}{\Gamma(1/6)}\int_0^\infty e^{-xs}s^{a-1}(1+s)^{c-a-1} ds$$
is a confluent hypergeometric function ($c > a >0$). We readily see that 
$g'(x) = - e^{-x}U_7(x)$ and $g''(x) = e^{-x}U_{10}(x),$ so that by
(\ref{LCE}) the MSU property for $X_{2/3}$ amounts to
$$(x^2 U_{10}(x) - x U_{7}(x)) U_4(x) \; \le\; x^2 U_7^2(x), \quad x\ge 0.$$ 
Using twice the contiguity relation (6.6.5) in \cite{E} and some simple
transformations, we then find the equivalence condition
\begin{equation}
\label{Whitt}
(xU_4(x) - U_1(x)/6)(U_7(x) - U_4(x))\;\ge\; -5U_4^2(x)/6, \quad x\ge 0.
\end{equation}
Notice that (\ref{Whitt}) is true as soon as $x\ge 1/6$ thanks to the obvious inequalities $U_7(x) \ge U_4(x)\ge U_1(x).$ However, some easy computations yield the asymptotics
$$U_7(x)\sim \frac{\Gamma(4/3)}{\Gamma(1/6)}x^{-4/3}, \;\; U_4(x)\sim
\frac{\Gamma(1/3)}{\Gamma(1/6)}x^{-1/3}\;\;\mbox{and}\;\; U_1(x) \to
\frac{\Gamma(2/3)}{\Gamma(5/6)}$$
 when $x\to 0^+,$ so that (\ref{Whitt}) does not hold anymore. This last discussion about the
case $\a = 2/3$ may look tedious, all the more that the proof that the MSU property does not hold for any $\a > 1/2$ is quite simple as we will soon see. But this must be read as a preparatory example for the sequel, since we will need inequalities such as (\ref{Whitt}) for some confluent hypergeometric functions in order to show the
MSU property when $\a\le 1/2,$ a fact which is more involved.

\section{Proof of the Theorem}

We begin with the only if part, and we will give three different arguments. The first one relies on the aforementioned Humbert-Pollard representation for $\fa$ - see e.g. formula (14.31) in \cite{S}:
$$f_\a (x)\; =\; \frac{1}{\pi} \sum_{n\ge 1} \frac{(-1)^{n-1}}{n!}\sin (\pi\a n) \Gamma(1+\a n) x^{-(1 +\a n)}\; =\; \sum_{n\ge 1} \frac{(-1)^n}{\Gamma (-n\a)n!} x^{-(1 +\a n)}.$$
Because $\a < 1$ this expansion may be differentiated term by term on $(0, + \infty),$  yielding
$$x f_\a' (x)\; =\;\sum_{n\ge 1} \frac{(-1)^{n-1}}{\Gamma (-n\a)n!}  (1+\a n)x^{-(1 +\a n)}$$
and
$$x^2 f_\a'' (x) \; +\; x f_\a' (x)\; =\; \sum_{n\ge 1} \frac{(-1)^{n-1}}{\Gamma (-n\a)n!}   (1+\a n)^2 x^{-(1 +\a n)}$$
for every $x > 0$ (these three series representations explain the
reformulation of (\ref{LCE}) in terms of a reverse Cauchy-Schwarz inequality
mentioned at the end of Section 1). Using the expansions up to $n = 2$ and the concatenation formula $\Gamma(z+1) = z\Gamma(z)$ we obtain, after some simplifications,
$$x^2(f_\a')^2(x)\, - \, (x^2 f_\a''(x) + xf_\a'(x))f_\a(x)\; =\; \frac{\a^2 x^{-(2 + 3\a)}}{2\Gamma(-\a)\Gamma(-2\a)}\; +\; {\rm o}(x^{-(2 + 3\a)})$$
in the neighbourhood of infinity, which entails that (\ref{LCE}) does not hold when $\a > 1/2,$ since the leading term is then negative. 

The second one hinges upon an expansion for the density $g_\a$ of the random
variable $Y_\a = \log \Za,$ which had been obtained in \cite{BB} - see (3.5)
and (6.4) therein - independently of the Humbert-Pollard formula:
$$g_\a(t)\; =\; e^{-\a t - e^{-\a t}}\sum_{j \ge 0} b_j \a^{j+1} (-1)^jR_j (-e^{-\a t})$$
for every $t\in\rl,$ where  the coefficients $b_j$ and $R_j(x)$ are defined through the entire series
$$\frac{1}{\Gamma(1+z)}\; =\; \sum_{j\ge 0} b_j z^j\quad \mbox{and}\quad e^{z + xe^z}\; =\; \sum_{j\ge 0} \lpa \frac{R_j(x) e^x}{j!}\rpa z^j.$$
Besides, setting 
$$P_\a(x) \; =\; \sum_{j \ge 0} b_j \a^{j+1} (-1)^jR_j (-x)$$
for any $x> 0,$ we know from (6.3) in \cite{BB} that the series converge absolutely, and with exactly the same argument one can show that it can be differentiated term-by-term. On the other hand, simple computations entail that $g_\a$  is log-concave over $\rl$ if and only if 
$$P_\a(x)^2 + x^2 (P'_\a(x))^2 \; \ge \; P_\a(x) (xP_\a'(x) + x^2P_\a''(x))$$
for every $x >0.$ Letting $x\to 0^+$ yields
$$P_\a(x)^2 + x^2 (P'_\a(x))^2\; \sim\; xP_\a(0)^2\; =\; x\lpa\frac{\a}{\Gamma(1-\a)}\rpa^2$$
and
$$P_\a(x) (xP_\a'(x) + x^2P_\a''(x))\; \sim\; xP_\a(0)P_\a'(0)\; =\; x\lpa\frac{\a}{\Gamma(1-\a)}\rpa^2\lpa 1 - \frac{\Gamma(1-\a)}{\Gamma(1-2\a)}\rpa$$
where the evaluations of $P_\a(0)$ and $P_\a'(0)$ rest upon the definition of $b_j$ and the fact that $R_j(0) = 1$ and $R_j'(0) = 2^j - 1.$ Similarly as above, we see that the second asymptotic is larger than the first one when $\a > 1/2.$

For the third and simplest argument, we will invoke an identity in law connecting two independent copies $Y_\a$ and ${\tilde Y}_\a$ of the random variable $\log \Za,$ which can be readily obtained in changing the variable in Exercise 4.21 (3) of \cite{CY}:
\begin{equation}
\label{ID}
Y_\a\, -\, {\tilde Y}_\a\; \elaw\; U_\a
\end{equation} 
where $U_\a$ is a real random variable with density
$$u_\a(x)\; =\; \frac{\sin \pi\a}{\pi(e^{\a x} + 2 \cos \pi\a + e^{-\a x})}\cdot$$
We compute then the second derivative of $x\mapsto \log (e^{\a x} + 2 \cos \pi\a + e^{-\a x})$ which is
$$\frac{4\a^2 + \a^2\cos \pi\a \cosh \a x}{(e^{\a x} + 2 \cos \pi\a + e^{-\a x})^2}\cdot$$
Hence, we see that $U_\a$ has a log-concave density over $\rl$ iff $\a \le 1/2.$ By Pr\'ekopa's theorem this entails that $\log \Za$ does not have a log-concave density over $\rl$ when $\a > 1/2,$ which means that $\Za$ is not MSU.

\begin{REMS} {\em (a) This negative result shows that the kernel $(x,y)\mapsto
  f_\a(e^{x-y})$ is not ${\rm TP}_2$ when $\a > 1/2,$ which contradicts the
  affirmation made in Lemma 1 (iv) of \cite{G} that this kernel is strictly totally
  positive for every $0<\a< 1$ (actually the contradiction could already have been seen
  in reading \cite{K} p. 390 carefully). Notice that this latter affirmation seems crucial to obtain the so-called 
  bell-shape property for all $\a$-stable variables with index $\a < 1$ - see
  p. 237 in \cite{G}. However, since this question is quite different from the topic of the present paper, we plan to tackle the problem (if really
  any) in some further research.

\vspace{2mm}

\noindent
(b) Though somewhat more technical, the methods resting upon
Humbert-Pollard's and Brockwell-Brown's expansions give some 
insight on the location where the inequality (\ref{LCE}) breaks down, an information which could not have been obtained by the third argument. I had also believed for a long time 
that these two expansions would give the if part, but this still eludes me because of the alternate signs.

\vspace{2mm}

\noindent
(c) It would be quite interesting to see if the third argument could not give the if part either. From the analytical viewpoint this would be the consequence of a positive answer to the following question. If $X$ is a real random variable with density such that the independent difference $X -X$ has a log-concave density, does $X$ have a log-concave density as well? This assertion, a kind of reverse to Pr\'ekopa's theorem for which we found neither references nor counterexamples in the literature, goes somehow the opposite way to the central limit theorem (which leads to a log-concave density after enough convolutions on any probability distribution with finite variance). Franck Barthe wrote me that he would be surprised if it held true in full generality. In our positive $\a-$stable framework, the fact
that $\log\Za$ is unimodal (and log-concave at both infinities when $\a\le 1/2$) might add some crucial properties, but overall I could not find any clue in the
direction of this statement.}
   
\end{REMS}

We now consider the if part, using yet another argument since the three
above methods turned out fruitless. We first prove two lemmas which are
of independent interest. The second one is a generalisation of Williams' result and could be formulated in several ways (we chose the one tailored to our purposes).

\begin{LEMO} Let $X$ be a ${\rm Beta}\, (\a, \b)$ variable and $Y$ an
  independent $\,{\rm Gamma}\, (c)$ variable such that $\b \le 1$ and $\a + \b \ge c.$ Then the product $X\times Y$ is MSU.
\end{LEMO}

\noindent
{\em Proof}: When $\a +\b = c,$ the result follows easily without further
assumption on $\beta$ because 
\begin{equation}
\label{IDG}
X\; \times\; Y \;\elaw\; {\rm Gamma} \, (\a),
\end{equation}
a fact which can be found e.g. in \cite{CY} Exercise 4.2. When $\a +\b \ge c,$
we first compute the density function of $X\times Y$: two changes of
variable entail
\begin{eqnarray*}
f_{X\times Y}(x) & = & \frac{\Gamma (\a +\b)}{\Gamma
  (\a)\Gamma(\b)\Gamma(c)}\int_0^1 e^{-x/u} (x/u)^{c-1} u^{\a -1} (1-u)^{\b
  -1} \frac{du}{u}\\
& = & \frac{\Gamma (\a +\b) x^{c-1}}{\Gamma
  (\a)\Gamma(\b)\Gamma(c)}\lpa e^{-x}\int_0^\infty e^{-xu} u^{\b-1}(u +1)^{c -(\a
  +\b)} du\rpa
\end{eqnarray*}
(notice that this recovers (\ref{IDG}) when $\a +\b = c$), so that we will be done as soon as the function $t\mapsto g_{\a,\b,c}(e^t)$ is
log-concave, with the notation
$$ g_{\a,\b,c}(x)\; =\;e^{-x}\int_0^\infty e^{-xu} u^{\b-1}(u +1)^{c -(\a
  +\b)} du.$$
Using now exactly the same discussion made at the end of Section 2 for the case $\a = 2/3$
(with adaptated computations) entails this log-concavity property is equivalent to
$$(x g_{\a,\b,c}(x) + (\a +\b -c)g_{\a,\b,c-1}(x)) (g_{\a,\b,c+1}(x) -
g_{\a,\b,c}(x))\; \ge \; (\b -1)(g_{\a,\b,c-1}(x))^2$$
for every $x\ge 0.$ But in the above, the right-hand side is negative because
$\b < 1,$ whereas the left-hand side is positive from the obvious inequality
$g_{\a,\b,c+1}(x)\ge g_{\a,\b,c}(x),$ and since by assumption $\a +\b \ge c.$    

\fin

\begin{LEMU} For all integers $p, n \ge 2$ such that $n > 2p,$ one has the following representation as an independent product:
$$Z_{p/n}^{-p}\; \elaw\; \frac{n^n}{p^p}\;{\rm Beta}\,(2/n, 1/p - 2/n)\,{\rm Gamma}(1/n)\;\times\;{\rm Beta}\,(4/n, 3/p - 4/n)\,{\rm Gamma}(3/n)\;\times\;\cdots\;$$

\vspace{-5mm}

\begin{flushright}

$\times\; {\rm Beta}\,(2(p-1)/n,
(p-1)/p - 2(p-1)/n)\,{\rm Gamma}((2p-3)/n)$

\vspace{3mm}

$\times\; {\rm
Gamma}\,((2p-1)/n) \;\times\;\cdots\;\times\;{\rm
Gamma}\,((n-1)/n).$
\end{flushright}
\end{LEMU}

\noindent
{\em Proof}: We first evaluate the fractional moments of $Z_{p/n}^{-p}$ using
(\ref{MellinS}), the duplication formula for the Gamma function - see
e.g. Formula (1.2.11) in \cite{E}, and some rearrangement involving the crucial assumption $n>2p:$ for every $s > - 1/n$ one obtains
\begin{eqnarray*}
\EE\lcr (Z_{p/n}^{-p})^s\rcr & = & \frac{\Gamma(ns +1)}{\Gamma(s+1)}\;
\times\;\frac{\Gamma(s +1)}{\Gamma(ps+1)} \\
& = & \lpa \frac{n^n}{p^p}\rpa^{\! s} \frac{\Gamma(s +1/n)\;\ldots\; \Gamma(s +
  (n-1)/n)\;\Gamma(1/p)\;\ldots\; \Gamma((p-1)/p)}{\Gamma(s +1/p\,)\;\ldots\; \Gamma(s +
  (p-1)/p)\;\Gamma(1/n)\;\ldots\; \Gamma((n-1)/n)}\\
& = & \lpa \frac{n^n}{p^p}\rpa^{\! s} \lpa\frac{\Gamma(s
  +2/n)\Gamma(1/p)}{\Gamma(s +1/p)\Gamma(2/n)}\rpa\lpa\frac{\Gamma(s
  +1/n)}{\Gamma(1/n)}\rpa
\end{eqnarray*}
\begin{eqnarray*}  
& & \;\;\;\;\;\;\;\;\;\;\;\;\;\;\;\;\;\;\;\;\;\;\;\;\;\;\;\;\;\times\;\lpa\frac{\Gamma(s
  +4/n)\Gamma(2/p)}{\Gamma(s +2/p)\Gamma(4/n)}\rpa\lpa\frac{\Gamma(s
  +3/n)}{\Gamma(3/n)}\rpa\;\times\;\cdots\\
& &
\;\;\;\;\;\;\;\;\;\;\;\;\;\;\;\;\;\;\;\;\;\;\;\;\;\;\;\;\;\;\;\;\;\;\;\;\times\;\lpa\frac{\Gamma(s
  + 2(p-1)/n)\Gamma((p-1)/p)}{\Gamma(s
  + (p-1)/p)\Gamma(2(p-1)/n)}\rpa\lpa \frac{\Gamma(s +
  (2p-3)/n)}{\Gamma((2p-3)/n)}\rpa\\
& & \;\;\;\;\;\;\;\;\;\;\;\;\;\;\;\;\;\;\;\;\;\;\;\;\;\;\;\;\;\;\;\;\;\;\;\;\;\;\;\;\;\;\;\times\;\lpa\frac{\Gamma(s + (2p-1)/n)}{\Gamma((2p-1)/n)}\rpa\;\times\;\cdots\;\times\;\lpa \frac{\Gamma(s +  (n-1)/n)}{\Gamma((n-1)/n)}\rpa.
\end{eqnarray*}

\noindent
On the other hand, it is well-known and easy to see that the fractional
moments of the ${\rm Beta}\, (\a, \b)$ and $\,{\rm Gamma}\, (c)$ variables are
given by
$$\EE\lcr ({\rm Beta}\, (\a, \b))^s\rcr\; =\; \frac{\Gamma(s+\a)\Gamma(\a
  +\b)}{\Gamma(s+\a +\b)\Gamma(\a)}\quad\mbox{and}\quad \EE\lcr ({\rm Gamma}\,
(c))^s\rcr\; =\; \frac{\Gamma(s+ c)}{\Gamma(c)}\cdot $$
The claim follows now by identification of the Mellin transform. 

\fin

\begin{REM} {\em As mentioned before, we see from this proof that analogous product representations for positive $\a$-stable laws with any $\a$ rational can be obtained accordingly. This might be useful to some other problems. Compare also with Theorem 2.8.4 in \cite{Z} where transforms of the densities $f_{p/q}$ are given as solutions to differential equations of higher order, an analytical representation which seems less tractable than Williams-type representations.}

\end{REM}

\noindent
{\bf End of the proof}: We need to show (\ref{LCE}) for any $\a\le 1/2$ and
$x \ge 0.$ Setting $\ga (x) = (x^2 f_\a''(x) + xf_\a'(x))f_\a(x) -
x^2(f_\a'(x))^2,$ we see from the Humbert-Pollard decomposition and its
differentiations that the application $\a\mapsto\ga(x)$ is continuous on
$(0,1)$ for every fixed $x \ge 0.$ By a density argument, it is therefore
sufficient to prove (\ref{LCE}) for any $\a = p/n$ with $p, n$ integers
greater than two such that $n> 2p,$ and every $x \ge 0.$ This amounts to the MSU
property for $Z_{p/n}^{-p},$ which now comes easily from Lemmas 1 \& 2, the MSU
property for Gamma variables and the stability of the MSU property by
independent multiplication.

\bigskip

\noindent
{\bf Acknowledgement.} This work was initiated during a stay at the University
of Tokyo. I am very grateful to Nakahiro Yoshida for his hospitality, and to
the grant ANR-09-BLAN-0084-01 for partial support.

\end{document}